\documentclass[12pt, a4paper, twosided, reqno]{amsart}
\usepackage{enumerate, cite}
\usepackage{hyperref}

\newtheorem{theorem}{Theorem}
\newtheorem{lemma}{Lemma}

\newtheorem{corollary}{Corollary}

\newtheorem{definition}{Definition}

\allowdisplaybreaks
\author[N. Gahlian]{Nidhi Gahlian }
\address{nidhi gahlian; department of mathematics, university of delhi, delhi-110007, india.}
\email{nidhigahlyan81@gmail.com}

\thanks {Research work of the author is supported by research fellowship from Department of Science and Technology(INSPIRE), New Delhi, India.}
\title[Oscillation Result]{Oscillation results of Higher order Linear differential equation}
\begin{document}
	\maketitle
	
\begin{abstract}
 We study higher order linear differential equation $y^{(k)}+A_1(z)y=0$ with $k\geq2$, where $A_1=A+h$, $A$ is a transcendental entire function of finite order with  $\frac{1}{2}\leq \mu(A)<1$ and $h\neq0$ is an entire function with $\rho(h)<\mu(A)$. Then it is shown that, if $f^{(k)}+A(z)f=0$ has a solution $f$ with $\lambda(f)<\mu(A)$ then exponent of convergence of zeros of any non trivial solutions of $y^{(k)}+A_1(z)y=0$ is infinite.
\end{abstract}
\section{\textbf{Introduction}}
For the understanding of this paper, we must know the basic facts of Nevanlinna's value distribution theory. For a meromorphic function $f$, $n(r,f)$, $N(r,f)$, $m(r,f)$ and  $T(r,f)$ denote un-integrated counting function, integrated counting function, proximity function and characteristic function respectively. We also use first main  theorem of Nevanlinna for a meromorphic function $f$, see \cite{hay,ilpo,yanglo}.
 Now we present elementary definitions of order of growth $\rho(f)$,  lower order of growth $\mu(f)$ exponent of convergence of zeros $\lambda(f)$ for a meromorphic function $f$ to make the paper self contained. 
$$\rho(f)=\limsup_{r\to\infty}\frac{\log T(r,f)}{\log r},$$
$$\mu(A)=\liminf_{r\to\infty}\frac{\log T(r,f)}{\log r},$$ 
 $$\lambda(f)=\limsup_{r\to\infty}\frac{\log n(r,\frac{1}{f})}{\log r}=\limsup_{r\to\infty}\frac{\log N(r,\frac{1}{f})}{\log r},$$ 
where $$N(r,\frac{1}{f})=\int_{0}^{r}\frac{n(t,\frac{1}{f})-n(0,\frac{1}{f})}{t}dt + n(0,\frac{1}{f})\log r.$$ 
A meromorphic function $g(z)$ is a small function of $f(z)$ if $T(r,g)=S(r,f)$ and vice versa, here $S(r,f)$ denote such quantities which are of growth $o(T(r,f))$ as $r\to\infty$, outside of a possible exceptional set of finite linear measure.\\
Suppose  $A$ be an entire function and $k\geq2$ is an integer. The complex differential  equation 
\begin{equation}\label{maineq1}
	f^{(k)} +A(z)f=0
\end{equation}
has entire solutions $f_j(j=1,2,...,k)$, as proved by Hille \cite{hile} that any solutions of \eqref{maineq1} are entire functions whenever $A$ is an entire function. 
In this note, we are concerned with the zero distribution of solutions of linear differential equation of $k^{th}$ degree with perturbed coefficients. The research in this direction is called complex oscillation theory and to classify the oscillations of solutions of \eqref{maineq1} has been a long standing problem since $1980s$. And it is basically depends on finding the conditions on $A(z)$ so that the existing solution $f$ of \eqref{maineq1} have $\lambda(f)=0$ or $\lambda(f)\geq\rho(A)$ or $\lambda(f)=\infty$.  However, rest cases like $0<\lambda(f)<\rho(A)$ or $\lambda(f)=\rho(A)=\infty$ or $\rho(A)<\lambda(f)<\infty$ are also possible but they seems to  be quite exceptional.\\
 In this field, many mathematicians have done research by relating the exponent of convergence of zeros $\lambda(f_j)(j=1,2,...,k)$ and  order of growth $\rho$ of coefficient $A$, for e.g \cite{alt,stebank,bank,lang,langly}. More results regarding complex oscillation of solutions of linear differential equation can be found in \cite{chen,chi,jklang} and references therein.

In $2005$, A. Alotaibi \cite{alt} proved the following theorem, which shows that by doing small perturbation of equation \eqref{maineq1} we get exponent of convergence of zeros of solution is  at least the order of growth of coefficient $A$ .
\begin{theorem}\label{theoremA}\cite{alt}
	Suppose that $A$ is a transcendental entire function with $\rho(A)<\frac{1}{2}$, $k\geq 2$ and 
\eqref{maineq1}
	has a solution $f$ with $\lambda(f)<\rho(A)$. Let 
	\begin{equation}\label{maineq2}
		A_1= A+h, 
	\end{equation} where $h\not\equiv 0$ is an entire function with $\rho(h)<\rho(A)$. Then exponent of convergence of zeros of any non trivial solution of
	\begin{equation}\label{eq3}
		g^{(k)}+A_1g=0
	\end{equation}
	does not have a solution $g$ with $\lambda(g)<\rho(A)$.
	\end{theorem}
 So, it seems interesting to find the condition on $A(z)$ so that exponent of convergence  of any nontrivial solution of equation \eqref{eq3} is infinite, and hence in $2020$, by using similar idea in \cite{alt}  J. Long and Y. Li proved the following theorem by considering lower order of growth of $A$ along with the small perturbation of such equation.
\begin{theorem}\label{theoremB}\cite{longli}
Suppose that $A$ is a transcendental entire function of finite order with $ \mu(A)<\frac{1}{2}$, $k\geq 2$ and 
\eqref{maineq1}
has a solution $f$ with $\lambda(f)<\mu(A)$. Let 
$A_1$ satisfies \eqref{maineq2}
 where $h\not\equiv 0$ is an entire function with $\rho(h)<\mu(A)$. Then exponent of convergence of zeros of any non trivial solution of \eqref{eq3} is infinite.
\end{theorem}
  Motivated by above results, It is natural to ask what we can say about the coefficient $A(z)$ of $\rho(A)\geq\frac{1}{2}$ or $\mu(A)\geq \frac{1}{2}$.
  So we consider the case $\frac{1}{2}\leq\mu(A)<1$ and prove the following result. 
\begin{theorem}\label{main1}
Suppose that $A$ is a transcendental entire function of finite order with $\frac{1}{2}\leq \mu(A)<1$, $k\geq 2$ and \eqref{maineq1}
 has a solution $f$ with $\lambda(f)<\mu(A)$. Let 
 $A_1$ satisfies\eqref{maineq2} and $h\not\equiv 0$ is an entire function with $\rho(h)<\mu(A)$. Then exponent of convergence of zeros of any non trivial solution of
\begin{equation}\label{maineq3}
 y^{(k)}+A_1y=0
\end{equation}
 is infinite.
\end{theorem}
By the proof of theorem 3 , we can easily prove the following result.
\begin{corollary}
	Let $A$ is a transcendental entire function of finite order with $\frac{1}{2}\leq \mu(A)<1$, $k\geq 2$ and \eqref{maineq1}
	has a solution $f$ with finitely many zeros, and let 
	$A_1$ satisfies\eqref{maineq2} and $h\not\equiv 0$ is an entire function with $\rho(h)<\mu(A)$. Then \eqref{maineq3} does not have a non trivial solution with finitely many zeros.
	\end{corollary}
 In section $2$ we state  some lemmas  and results. In section 3, we prove Theorem $3$.
\section{\textbf{Auxiliary results}}
In this section we present some lemmas and definition which will be helpful in proving our main theorem. For a set $I\subset (0,\infty)$, the linear measure is defined by $m(I)=\int_{I}\,dt$. For a set $J\subset (1,\infty)$, the logarithmic measure is defined by $m_{l}(J)=\int_{J}\frac{1}{t}\,dt$ and the upper and lower logarithmic density of $J\subset[1,\infty)$ by
$$\overline{logdens}J= \limsup_{r\to\infty}\frac{\ m_l(J\cap[1,r])}{\ log}$$ and $$\underline{logdens}J= \liminf_{r\to\infty}\frac{\  m_l(J\cap[1,r])}{\ log}.$$ The logarithmic density actually gives an idea how big the set J is.

\begin{definition} \cite{ilpo} 
Let $B(z_n,r_n)= \{ z: |z- z_n|<r_n\}$ be the open disc in the complex plane.Countable union $\bigcup\limits_{n=1}^{\infty} B(z_n,r_n)$ is said to be $R$-set if $z_n\to\infty$ and $\Sigma r_n$ is finite.

\end{definition}
Now we state a very well known result regarding rational function, which we use  in our main theorem's proof.
\begin{lemma}\label{impleA}\cite{ilpo}
	A meromorphic function $f$ is a rational function iff  $T(r,f)=O(\log r)$.
\end{lemma}
Next we define a well known representation for higher order logarithmic derivatives
through which we can easily show the possibility of the existence of a solution $f$ of \eqref{maineq1} with no zeroes by taking $F=\frac{\ f'}{\ f}, f=e^P$ where $P$ is an entire function.
\begin{lemma}\label{imple1}(Hayman's Lemma)\cite{hay}
Let $f(z)$ be an analytic function, and let $F= \frac{\ f'}{\ f}$. Then for $k\in \mathbb{N}$, we have 
$$\frac{\ f^{(k)}}{\ f}= F^k + \frac{\ k(k-1)}{\ 2}F^{k-2}F'+ P_{k-2}(F),$$ where $P_{k-2}$  is a differential polynomial with constant coefficients, which vanishes identically for $k\leq2$ and has degree of $(k-2)$ when $k>2$.
\end{lemma}
\begin{lemma}\label{imple2}\cite{yang}
	Let $B(z)$ is an entire function with $\mu(B)\in \left[\frac{1}{2},\infty\right)$ then there exists a sector $\Omega(\alpha, \beta),$ $\beta -\alpha\geq \frac{\pi}{\mu(B)}$, such that $$\overline{ \lim_{r\to\infty}} \frac{\log\log|B(re^{\iota\theta})|}{ \log r} \geq \mu(B)$$ $\forall \theta\in \Omega(\alpha,\beta)$, where $0\leq\alpha<\beta<2\pi.$
\end{lemma}
\begin{lemma}\label{imple3}\cite{ilpo}
	Suppose that $f(z)$ is a meromorphic function of finite order. Then there exists a positive integer $N$ such that 
	$$\frac{f'(z)}{f(z)}=O(|z|^N)$$
	holds for large $z$ outside of an $R$-set. 
	\end{lemma}


\begin{lemma}\label{imple4}(Langley's theorem)\cite{lang}
 Let $A(z)$ be a transcendental entire function of finite order, and let $J_1$ be a subset of $[1,\infty)$ of infinite logarithmic measure and with the following property. For each $r\in J_1$, there exists an arc $$ a_r= \{ re^{it} : 0\leq\alpha_r\leq t\leq\beta_r\leq 2\pi \}$$ of the circle $S(0,r)=\{ z:|z|=r\}$ such that 
	$$\lim_{r\to\infty,r\in J_1} \frac{\min \{\log|A(z)|: z\in a_r \}}{ \log r}= + \infty.$$
	Let $k\geq 2$ and let $f$ be a solution of \ref{maineq1} with $\lambda(f)<\infty.$ Then there exists a subset $J_2\subset[1,\infty)$ of finite measure, such that for large $r\in J_0 =J_1\setminus J_2,$  we have $$\frac{\ f'}{\ f}=c_rA(z)^{1/k} - \frac{\ k-1}{\ 2k}\frac{\ A'(z)}{\ A(z)} +O(r^{-2})$$ holds for all $z\in a_r$, where the constant $c_r$ satisfies $c_r^k=-1$ and may depend on r, for a given $r\in J_0$ but not depend on $z$, and the branch of $A(z)^{1/k} $ is analytic on $a_r$(included in the case where $a_r$ is the whole circle $S(0,r))$.
	
\end{lemma}

\section{\textbf{Proof of main theorem}}

\begin{proof}
Given $A_1= A+h$, $\rho(h)<\mu(A)\leq\rho(A)$ implies $\rho(A_1)=\rho(A)$. Let equation \eqref{maineq1} has a solution $f$ with $\lambda(f)<\mu(A)\leq \rho(A)$ and let us assume \eqref{maineq3} has a solution $y$  with $\lambda(y)<\infty.$ So we can take $$f= Pe^U$$ and $$ y=Qe^V,$$ where $U, V, P,Q $  are entire functions of finite order\cite{jklang}.\\
	 Now $f= Pe^U$ implies $\lambda(f)= \rho(P)<\infty $\cite{ilpo} as $e^U\neq 0$ \& $\rho(P)=\lambda(P).$ \\Similarly, $$\lambda(y)=\rho(Q).$$
Let 
\begin{equation}\label{eq4}
F=\frac{f'}{f}, 
\end{equation}
\begin{equation}\label{eq5}
Y= \frac{y'}{y}.
\end{equation}
	\\Using $f=Pe^U$ $\&$ $Y=Qe^V$, we get 
	
	$$ F= \frac{\ Pe^{U}U'+e^UP'}{\ Pe^U}$$
	\begin{equation}\label{eq6}
	=\frac{P'}{P} +U'. 
	\end{equation}
	Similarly,
	\begin{equation}\label{eq7}
	 Y=\frac{Q'}{Q}+V'.
	\end{equation}
	Applying Hayman's Lemma \ref{imple1}, we get 
	\begin{equation}\label{eq8}
\frac{f^{(k)}}{f}=F^k + \frac{k(k-1)}{2} F^{k-2}F'+ P_{k-2}(F)
	\end{equation} 
and
\begin{equation}\label{eq9}
	 \frac{y^{(k)}}{y}= Y^k +\frac{k(k-1)}{2}Y^{k-2}Y'+ P_{k-2}(Y),
\end{equation}
 where $P_{k-2}$  is a differential polynomial with constant coefficients, which vanishes identically for $k\leq2$ and has degree of $k-2$, when $k>2$.\\
 Choose
\begin{equation}\label{eq10}
		max \{\lambda(f),\lambda(y), \rho(h)\}< \beta<\gamma<\mu<1. 
\end{equation}
	\\
	Using lemma \ref{imple2} for $A(z)$, there exists a sector 
	$\Omega(\alpha,\beta)$
	such that following inequality holds $\forall \theta\in \Omega(\alpha,\beta)$.

	$$\overline {\lim_{r\to\infty}}\frac{\log\log|A(re^{\iota\theta})|}{\log r}\geq\mu(A).$$ 
	This gives 
	$$\log|A(re^{\iota\theta})|\geq r^\gamma,$$ where $\gamma=\mu(A)-\epsilon$.
Set	$$ J_1:=\{z=re^{\iota\theta}:|z|=r>r_0, \alpha_0 < \theta< \beta_0\},$$ where $\alpha<\alpha_0<\beta_0<\beta$ and $r_0$ is a fixed number, satisfying
\begin{equation}\label{eq11}
	\inf_{|z|=r\in J_1 }\log |A(z)|\geq r^\gamma,
\end{equation}
where $J_1$ has a positive upper logarithmic density\cite{mani}.
	 \\Using lemma \ref{imple3}, there exist a set $J_2\subset[1,\infty)$ which is a subset of finite measure, for some $t_1\in\mathbb{N},$

\begin{equation}\label{eq12}
\left| \frac{A'(z)}{A(z)} \right| +\left| \frac{P'(z)}{P(z)}\right|+ \left|\frac{Q'(z)}{Q(z)}\right| \leq r^\tau, 
\end{equation}
holds for $|z|=r\geq1 ,  r\notin J_2$.	 Now for large $r\in J_1$ and $\rho(h)<\rho(A)$, we have
	 $A_1=A+h$. This gives
\begin{align*}
\log(A_1)&=\log(A+h)\\
& =\log A +o(1)\\
& \geq r^\gamma +o(1).     
\end{align*}

Next we calculate $\frac{f'}{f}$ and $\frac{y'}{y}$ in terms of $A(z)$.
 For applying lemma \ref{imple4}, take an arc  $a_r := \{ z=r_1e^{\iota\theta} :\theta\in(\alpha_1,\beta_1)\}$ in $J_1$ for some fixed $r_1$, where $\alpha_0<\alpha_1<\beta_1<\beta_0$. Now for some   $\theta_1 \in (\alpha_1,\beta_1)$, we get $$\min \{\log |A(z)|:z\in a_r\}=\log |A(re^{\iota\theta_1})|.$$
 Next, $$ \frac{\log A(re^{\iota\theta_1})}{\log r}\geq\frac{r^\gamma}{\log r}.$$
 As $\gamma<\mu<1$, so $$ \lim_{r\to\infty}\frac{r^\gamma}{\log r} \rightarrow\infty.$$ 
 Hence $$\lim _{r\to\infty}\frac{\ min \{ \log |A(z)|:z\in a_r\} }{\log r}=\infty.$$
Next, on applying lemma \ref{imple4} in equations \eqref{maineq1} and  \eqref{maineq3}. The following equalities hold for large $r\in J_0$. 

 \begin{equation}\label{eq13}
 	\frac{f'}{f}= cA(z)^{1/k} -\frac{k-1}{2k} \frac{A'(z)}{A(z)}+ O(r^{-2}),  z\in \Omega, c^k=-1,
 \end{equation}
 \begin{equation}\label{eq14}
 	\frac{y'}{y}= dA_1(z)^{1/k}-\frac{k-1}{2k}\frac{A'_1(z)}{A(z)}+O(r^{-2}),    z\in \Omega, d^k=-1.
 \end{equation}
 Now expanding $A_1(z)^{1/k}$ and $\frac{A'_1(z)}{A(z)}$ in terms of $A(z)^{1/k}$ and $\frac{A'(z)}{A(z)}$ with the help of binomial theorem, we have
 
 \begin{align*}
 &\limsup_{r\to\infty}\frac{\log\log M(r,h)}{\log r}=\rho(h),\\
 \implies
 &|h(z)|\leq e^{r^{\rho(h)+o(1)}},
\end{align*}
where $M(r,h) $ is the maximum term.

 For $|z|=r\to\infty ,r\in J_0$, using above equation and equation \eqref{eq11}, we get $$\left|\frac{h(z)}{A(z)}\right|\leq \frac{e^{r^{\rho(h)+0(1)}}}{e^{r^\gamma}} = o(1).$$ Similarly we can find $$\left|\frac{h'(z)}{A(z)}\right|\leq\frac{e^{r^{\rho(h)+0(1)}}}{e^{r^\gamma}} = o(1).$$
 Now expanding with the help of above two equations, we have
\begin{align*} 
A_1^{\frac{1}{k}}(z)&=(A+h)^{\frac{1}{k}}\\
 &=A^{\frac{1}{k}}\left(1+ \frac{h}{A}\right)^{\frac{1}{k}}\\
 &= A^{\frac{1}{k}}\left(1+ O\left({\frac{|h|}{|A|}}\right)\right),
\end{align*}
for $|z|=r\in J_0$. Similarly,
\begin{align*}
 \frac{A_1'(z)}{A(z)}&=\left(\frac{A'+h}{A+h}\right)\\&=\frac{A'+h}{A\left(1+\frac{h}{A}\right)}\\&=\frac{A'+h'}{A}\left(1-\frac{h}{A}
+\frac{h^2}{A^2}...\right)\\&=\left(\frac{A'}{A}+\frac{h'}{A}\right)\left(1+O\left( \frac{|h|}{|A|}\right)\right)\\&=\frac{A'}{A} \left(1+O\left(\frac{|h|}{|A|}\right)\right),
\end{align*}
for $|z|=r\in J_0$.
So we get equation
\begin{equation}\label{eq15}
\frac{y'}{y}=dA(z)^{1/k}-\frac{k-1}{2k}\frac{A'(z)}{A(z)}+O(r^{-2}), d^k=-1.
\end{equation}

Now we will prove $c=d$ for large $r\in J_0.$
Let $d/c=\omega$, where $\omega^k=1$ and \\using equation \ref{eq10}, we get
\begin{equation}\label{eq16}
\frac{y'}{y}= c\omega A(z)^{1/k}-\frac{k-1}{2k}\frac{A'(z)}{A(z)}+O(r^{-2}), w^k=1.
\end{equation}
 Multiply equation \ref{eq8}  with $\omega$,
 $$\omega\frac{f'}{f}=\omega c A(z)^{1/k}-\omega \frac{k-1}{2k}\frac{A'(z)}{A(z)}+O(r^{-2}), c^k=-1.$$
Now subtract equation \ref{eq16} from above equation, we get
$$\omega \frac{f'}{f}-\frac{y'}{y}=\omega cA(z)^{1/k}-\omega \frac{k-1}{2k}\frac{A'(z)}{A(z)}+O(r^{-2})-\omega c (A(z)^{1/k})-\frac{k-1}{2k}\frac{A'(z)}{A(z)}+O(r^{-2})$$
$$\omega \left(\frac{f'}{f}-\frac{k-1}{2k}\frac{A'(z)}{A(z)}\right)=\frac{y'}{y}+ \frac{k-1}{2k}\frac{A'(z)}{A(z)}+O(r^{-2})$$
Now  by using Argument principle, integrate above equation around $|z_n|=r_n\in J_0,$ and $r_n\to\infty$ as $n\to\infty$, we get
\begin{equation}\label{eqA}
\omega\left[2\pi\iota   n\left(r_n,\frac {1}{f}\right) +\frac{k-1}{2k}2\pi\iota n\left(r_n,\frac{1}{A}\right)\right]+o(1)=2\pi\iota n\left(r_n,\frac{1}{y}\right)+\frac{k-1}{2k} 2\pi\iota n\left(r_n,\frac{1}{A}\right).
\end{equation}
 In this equation R.H.S must  be a positive integer as $ n\left(r_n,\frac{1}{y}\right)\geq 0$ and $n\left(r_n,\frac{1}{A}\right)>0.$
 Let us suppose that $n\left(r_n,\frac{1}{A}\right)=0$,  if so then it implies $N\left(r_n,\frac{1}{A}\right)=0$.\\ Now since $ \log|A(z)|\geq r^{\gamma}$ for $r>r_0$ i.e  
 $\inf_{|z|=r_n\in J_0 } \log |A(z)|$ is very large for $r_n\to\infty$, thus we get
\begin{align*}
 m(r_n,\frac{1}{A})&=\frac{1}{2\pi}\int_{0}^{2\pi}\log\left|\frac{1}{A(re^{\iota\theta})}\right|d\theta\\&=0.
 \end{align*}
 Hence $$T(r_n,\frac{1}{A})=m(r_n,\frac{1}{A})+n(r_n,\frac{1}{A})=0.$$
Now with the help of  Nevanlinna's first fundamental theorem 
$ T(r_n,A)=T(r_n,\frac{1}{A})+O(1)$, we get $$T(r_n,A)=O(1),$$
which is a contradiction  with the fact that A is transcendental. So our supposition is wrong and hence $n(r_n.\frac{1}{A})>0$.
So $n(r_n,\frac{1}{A})+n(r_n,\frac{1}{f}) $ is a positive integer.
Now as $n(r_n,\frac{1}{A}) \geq1$ implies $ \frac{k-1}{2k}n(r_n,\frac{1}{A})+n(r_n,\frac{1}{f}) \geq \frac{k-1}{2k}. $ 
Taking modulus of imaginary part of both sides after multiplying with $\omega$, we get
\begin{equation*}
	\left| Im\left[\omega \left( n\left (r_n,\frac{1}{f}\right)+\frac{k-1}{2k} n\left(r_n,\frac{1}{A}\right)\right)\right]\right|\geq \frac{k-1}{2k}|Im(\omega)|,
	\end{equation*}
 and taking imaginary part of both sides of  $\ref{eqA}$, we get
 \begin{equation*}
 	Im\left[\omega \left( n\left (r_n,\frac{1}{f}\right)+\frac{k-1}{2k} n\left(r_n,\frac{1}{A}\right)\right)\right] +Im(o(1)) =0,
\end{equation*}
so 
\begin{align*}
	|Imo(1)|=& \left|-Im\left[\omega \left( n\left (r_n,\frac{1}{f}\right)+\frac{k-1}{2k} n\left(r_n,\frac{1}{A}\right)\right)\right]\right|\\&
	\geq\frac{k-1}{2k}|Im(\omega)|\\& \geq \frac{k-1}{2k}\Delta,
	\end{align*}
 where $\Delta:= inf\{|Im(\omega)|:\omega^k=1,Im(\omega)\neq0\}.$ 
 It is obvious $\Delta >0.$ 
 Now for sufficiently large $r_n\in J_0$ and $|Im(o(1))|<\Delta \frac{k-1}{2k}$, 
 we get $Im(\omega)=0$.
\\ As $\omega^k=1$ with imaginary part zero implies either $\omega =-1$ or $\omega=1$.\\
 But $\omega=-1$  contradicts with  \eqref{eqA}, so $\omega =1$ and hence $c=d$.
 By using it we can write equation \eqref{eq15} as
\begin{equation}\label{eq17}
Y(z)=\frac{y'}{y}=cA(z)^{1/k}-\frac{k-1}{2k}\frac{A'(z)}{A(z)}+O(r^{-2}), c^k=-1,
\end{equation}
for $|z|= r\in J_0$.
Subtracting equation \eqref{eq16} from  equation \eqref{eq13}, we get 
$$\frac{f'(z)}{f(z)}=\frac{y'(z)}{y(z)}+O(1), |z|=r.$$
Now as  poles of $\frac{f'}{f}=$ zeroes of $f$, hence for large $r\in J_0$, we have
$$n(r,\frac{1}{f})= n(r,\frac{1}{y}).$$
 With the help of equations \eqref{eq5}, \eqref{eq6}, \eqref{eq7} and \eqref{eq8} we get,
$$\frac{P'}{P}+ U'=\frac{Q'}{Q}+V'+O(1).$$
Using  above equation and \eqref{eq12}, we get
$$|U'-V'|\leq \left|\frac{P'}{P}\right|+\left|\frac{Q'}{Q}\right|$$
$$\leq 2r^\tau,$$  which holds for all $|z|=r $ with large $ r\in J_0\setminus J_2$. Next 
$$\log M(r,U'-V')\leq \log(2r^{\tau}),$$ where $M(r,U'-V')$ is the maximum term, this gives
$$\frac{T(r,U'-V')}{\log r}\leq \tau.$$
As  $U'-V'$  is entire function so $$M(r,U'-V')=T(r,U'-V').$$\\
By above inequality, we get
$$T(r,U'-V')=O(\log r).$$ This implies $U'-V'$ is a rational function but 
as $U$ and $V$ are entire so $P_0 =U'-V'$ is a polynomial. From equation \eqref{eq7} and \eqref{eq8}, $$F=\frac{P'}{P}+U'=\frac{P'}{P }+P_0+V'$$ and
$$Y=\frac{Q'}{Q}+V'.$$  This implies
\begin{equation}\label{eq18}
	F=\frac{P'}{P}+P_0+Y-\frac{Q'}{Q}=Y+M,
\end{equation}
where $M= \frac{P'}{P}-\frac{Q'}{Q}+P_0.$
Using \eqref{maineq1} and \eqref{eq9}, we get
\begin{equation}\label{eq19}
 \frac{f^{(k)}}{f}=F^k + \frac{k(k-1)}{2} F^{k-2}F'+ P_{k-2}(F)=-A,
\end{equation}
 where $P_{k-2}$  is a differential polynomial with constant coefficients, which vanishes identically for $k\leq2$ and has degree of atmost $k-2$, when $k>2$.\\
By equations \eqref{maineq2}, \eqref{maineq3} and \eqref{eq10}, we get
\begin{equation}\label{eq20}
 \frac{y^{(k)}}{y}= Y^k +\frac{k(k-1)}{2}Y^{k-2}Y'+ P_{k-2}(G)=-A-h.
\end{equation}
Using  equations \eqref{eq18} and \eqref{eq19}, we get
$$(Y+M)^k+\frac{k(k-1)}{2}(Y+M)^{k-2}(Y'+M')+P_{k-2}(Y+M)=-A.$$
With the help of Binomial theorem, we can expand above equality and we obtain
$$Y^k+MkY^{k-1}+\frac{k(k-1)}{2}Y^{k-2}Y'+ B_{k-2}(Y,M)=-A,$$
where $B_{k-2}$ represent a polynomial in M,Y, and their derivatives with the total degree of atmost $k-2$.
Now combining above equation with  equation \eqref{eq20}, we get
\begin{equation}\label{eq21}
h= kMY^{k-1} +R_{k-2}(Y,M).
\end{equation}
Now we will claim  $M\not\equiv 0$.
On the contrary, let $M\equiv0$. As $F=Y+M$ implies $F=Y$, then by equations \eqref{eq19} and \eqref{eq20}, we get $h=0$, which contradicts with the hypothesis and hence claim is true.\\
Now divide equation \eqref{eq21} by $MY^{k-2}$, we get
\begin{equation}\label{eq22}
kY+\frac{R_{k-2}(Y,M)}{MY^{k-2}}=\frac{h}{MY^{k-2}}.
\end{equation}
Assume that  $|Y|>1$ and  $\frac{R_{k-2}(Y,M)}{MY^{k-2}}$ is a sum of the terms
$$\frac{1}{MY^{k-2}}M^{p_0}(M')^{p_1}...(M^{(k)})^{p_k}Y^{q_0}(Y')^{q_1}...(Y^{(k)})^{q_k},$$
where $q_0+q_1+...q_k\leq k-2$.\\
 As $|Y|>1$ implies $\frac{1}{|Y|}<1$ and taking modulus of  above equation, we get
$$|M|^{p_0+p_1+...p_k-1}\left|\frac{M'}{M}\right|^{p_1}.......\left|\frac{M^{(k)}}{M}\right|^{p_k}|Y|^{q_0+q_1+...q_k-k+2}\left|\frac{Y'}{Y}\right|^{q_1}...\left|\frac{	Y^{(k)
}}{Y}\right|^{q_k}$$ 

\begin{equation}\label{eq23}
\leq|M|^{p_0+p_1+...p_k-1}\left|\frac{M'}{M}\right|^{p_1}...\left|\frac{M^{(k)}}{M}\right|^{p_k}\left|\frac{Y'}{Y}\right|^{q_1}...\left|\frac{	Y^{(k)}}{Y}\right|^{q_k}.
\end{equation}
Taking proximity function on both sides of  equation \ref{eq22}, we get
\begin{align*}
m(r,kY)= m\left(r,-\frac{R_{k-2}}{MY^{k-2}} + \frac{h}{MY^{k-2}}\right)\\ m(r,Y)\leq m\left(r,-\frac{R_{k-2}}{MY^{k-2}} \right) +m\left( r,\frac{h}{MY^{k-2}}\right) +\log 2.\\&
\end{align*}
Using  equation \eqref{eq23}, we get
\begin{align*}
	m(r,Y)\leq& m\left( r,|M|^{p_0+ p_1+...p_k-1}\right)+m\left(r,\left|\frac{M'}{M}\right|\right)+...m\left(r,\left|\frac{M^{(k)}}{M}\right|\right)\\&+
	m\left(r,\left|\frac{Y'}{Y}\right|\right)...m\left(r,\left|\frac{Y^{(k)}}{Y}\right|\right)+ m(r,h)+m\left(r,\frac{1}{M}\right)\\&+m\left(r,\frac{1}{Y^{k-2}}\right)+\log4.
\end{align*}
  Now as M and Y are rational functions, so $m\left(r,\frac{M'}{M}\right) = S(r,M)$and $m\left(r,\frac{Y'}{Y}\right) = S(r,Y)$, hence we get
\begin{align*}
	m(,Y)\leq & c_0m(r,M)+m(r,\frac{1}{M})+m(r,h)+s(r,M)+S(r,Y),
\end{align*}
where $c_0=p_0+p_1+...p_k-1$, positive constant. Now by adding and subtracting terms  $c_0N(r,M),$ $N(r,\frac{1}{M}),N(r,h)$ and using first fundamental theorem of Nevanlinna, we get
\begin{align*}
	m(r,Y)\leq& c_0T(r,M)+T(r,\frac{1}{M})+T(r,h)+S(r,Y)
\end{align*}
\begin{align}\label{eq24}
	=(c_0+1)T(r,M)+T(r,h)+S(r,Y).
\end{align}\\ 
Take proximity  function on both sides of equation \eqref{eq20}, we get 
\begin{align*}
m(r,A)\leq& m(r,Y^k) +m(r,\frac{k(k-1)}{2} Y^{k-2}Y')+ m(r,P_{k-2}(Y)) +m(r,h)\\ &
\leq km(r,Y)+m(r,\frac{k(k-1)}{2})+m(r,Y^{k-2})+m(r,\frac{Y'}{Y})\\&+ m(r,Y)+S(r,Y)+ S(r,h).
\end{align*}
With some simple calculations, we get
\begin{align}\label{eqN}
	m(r,A)\leq c_2m(r,Y)+O(\log r).
	\end{align}
Using equation \eqref{eq24} and \eqref{eqN}, we get
\begin{align}\label{eq25}
	T(r,A)\leq c_3T(r,M)+ c_2T(r,h)+ S(r,Y) + O(\log r).
\end{align}
As
\begin{align*}
T(r,M)&= T\left(r,\frac{P'}{P}-\frac{Q'}{Q}+P_0\right)\\ &
=m\left(r,\frac{P'}{P}-\frac{Q'}{Q}+P_0\right)+N\left(r,\frac{P'}{P}-\frac{Q'}{Q}+P_0\right)\\ &
\leq S(r,P)+ S(r,Q)+O(\log r)+N\left(r,\frac{P'}{P}\right)+ N\left(r,\frac{Q'}{Q}\right)+ N(r,P_0)\\&
=N(r,\frac{1}{P})+N(r,\frac{1}{Q})+ S(r,W)\\&\leq T(r,P)+T(r,Q)+S(r,W),
\end{align*}

where W is some entire function with order of growth $\beta$.
Using equation \eqref{eq11} with $\rho(P)=\lambda(f)$ and $\rho(Q)=\lambda(y)$, we have
$$T(r,M)\leq o(r^{\beta}).$$
Now as $\rho(h)<\beta$, we get
$$T(r,h)\leq o(r^{\beta}).$$
Using these in equation \eqref{eq25}, we get 
$$T(r,A)\leq o(r^{\beta})$$  for $r\in J_0 \setminus J_2.$
Hence contradiction arises. This completes the proof.
\end{proof} 

\end{document}